\documentclass[a4paper,11pt]{article}
\usepackage{amsmath,amssymb,amsthm,a4wide,color}
\usepackage[utf8]{inputenc}
\usepackage[T1]{fontenc}
\usepackage{authblk,mathrsfs}
\usepackage{algorithm,algpseudocode}
\usepackage{graphicx,caption,subcaption}
\usepackage{textcomp}
\usepackage{}

\usepackage{booktabs}

\usepackage{caption}
\usepackage{subcaption}

\usepackage{pdfsync}
\newcommand{\RR}{\mathbb{R}}
\newcommand{\CC}{\mathbb{C}}
\newcommand{\mrm}{\mathrm}
\newcommand{\mL}{\mathrm{L}}
\newcommand{\mH}{\mathrm{H}}
\newcommand{\mJ}{\mathrm{J}}

\newcommand{\Id}{\mathrm{Id}}

\newcommand{\mT}{\mathrm{T}}
\newcommand{\mS}{\mathrm{S}}

\newcommand{\mR}{\mathrm{R}}

\newcommand{\Vh}{\mathrm{V}_{h}}

\newcommand{\mbVh}{\mathbb{V}_{h}}
\newcommand{\mbXh}{\mathbb{X}_{h}}

\newcommand{\bp}{\boldsymbol{p}}

\newcommand{\bx}{\boldsymbol{x}}
\newcommand{\bb}{\boldsymbol{b}}
\newcommand{\br}{\boldsymbol{r}}
\newcommand{\bd}{\boldsymbol{d}}

\newcommand{\bn}{\boldsymbol{n}}
\newcommand{\bu}{\boldsymbol{u}}
\newcommand{\bv}{\boldsymbol{v}}
\newcommand{\bw}{\boldsymbol{w}}
\newcommand{\bq}{\boldsymbol{q}}
\newcommand{\bl}{\boldsymbol{l}}
\newcommand{\bg}{\boldsymbol{g}}

\newcommand{\Bh}{\mathrm{B}}
\newcommand{\Ah}{\mathrm{A}}

\newcommand{\loc}{\mathrm{loc}}
\newcommand{\Prec}{\mrm{P}}

\newtheorem{defn}{Definition}[section]
\newtheorem{lem}[defn]{Lemma}

\newtheorem{prop}[defn]{Proposition}

\title{\textbf{Accelerating non-local exchange in}\\
  \textbf{generalized optimized Schwarz methods}}
\author[1]{X.Claeys}
\author[1,2]{R.Atchekzai}

\affil[1]{Sorbonne Université,
  Laboratoire Jacques-Louis Lions}

\affil[2]{CEA, CESTA}

\begin{document}

\maketitle

\begin{abstract}
  The generalized optimised Schwarz method proposed in [Claeys \&
    Parolin, 2022] is a variant of the Després algorithm for solving
  harmonic wave problems where transmission condition are enforced by
  means of a non-local exchange operator. We introduce and analyse an acceleration technique that significantly
  reduces the cost of applying this exchange operator without deteriorating the precision
  and convergence speed of the overall domain decomposition algorithm.
\end{abstract}

\section*{Introduction}

The present article is concerned with the efficient numerical solution
of time harmonic scalar wave equation by means of a substructuring
domain decomposition method. In a recent series of contributions, we
introduced variants of Després' algorithm, also dubbed Optimized
Schwarz Method (OSM), able to treat the presence of cross points in a
systematic manner while maintaining geometric convergence of the
overall DDM algorithms. This approach was proved a generalization of
classical OSM in the sense that it coincides with it
under appropriate circumstances. A full convergence framework was
also provided for this new approach, including a precise
quantification of convergence rates.

In classical OSM, wave equations are solved locally in each subdomain.
The local solves are then coupled by means of a swapping operator
$\Pi_{\loc}$ that exchanges ingoing/outgoing traces through each
interface. In the generalized variant of OSM introduced in
\cite{claeys2020robust}, a key innovation lies in a more sophisticated
exchange operator $\Pi$ that replaces the swapping operator. While
$\Pi_{\loc}$ is local by nature, the new exchange operator $\Pi$ is
non-local because it a priori couples distant non-neighbouring
subdomains (although $\Pi=\Pi_{\loc}$  in well identified circumstances).

Compared to the standard OSM, the generalized OSM leads to fastly
converging DDM algorithms, but requires dealing with a
potentially non-local exchange operator instead of the initial
swapping operator, which represents an extra non-negligible
computational cost. To be more precise, while the operation
$\bx\to\Pi_{\loc}(\bx)$ is trivial and simply consists in a
permutation of unknowns, the operation $\bx\to \Pi(\bx)$ requires the
solution to a global problem that has nevertheless the favorable
property of being hermitian positive definite. The goal of the present
article is to exhibit one strategy that allows to perform the exchange
operation $\bx\to \Pi(\bx)$ approximately but much faster.  We will
prove in addition that this approximation does not induce any error in
the overall DDM algorithm.

\quad\\
The exchange operation $\bx\to \Pi(\bx)$ requires solving a
self-adjoint positive definite (SPD) linear system, for which we rely
on a preconditioned conjugate gradient (PCG) solver.  At each
iteration $n$ of the DDM algorithm, such a linear system
$\mrm{L}(\bx^{(n)}) = \bb^{(n)}$ has to be solved. Two remarks can be
made that allow to substantially accelerate these linear solves.
First of all, the linear operator $\mL$ is independent of $n$. Besides, the
right hand sides $\bb^{(n)}$ vary from one step $n$ to another, but they
form a converging sequence because of the convergence of the overall DDM
algorithm. In this context, our acceleration strategy
consists in a simple recycling strategy combined with a brutal
truncation of PCG (only a few iterations are needed). After a
description of our method, we shall give numerical evidence of the
performance of this approach. In the last of this contribution,
we give a theoretical justification of the efficiency through
derivation of an explicit convergence. 

\section{Scattering problem under study}
We start by describing a typical wave propagation boundary value problem.
The aim of the domain decomposition we are discussing here is to solve
this problem as efficiently as possible. In the sequel $\Omega\subset \RR^d$
will refer to a polygonal/polyhedral bounded domain, and we wish to solve the
boundary value problem 
\begin{equation}\label{InitPb}
  \begin{aligned}
    &\text{Find}\;u\in\mH^{1}(\Omega)\;\text{such that}\\
    &\Delta u + \kappa^{2}u = 0\;\;\text{in}\;\Omega,\\
    &\partial_{\bn}u -i\kappa u= f\;\;\text{on}\;\partial\Omega. 
  \end{aligned}
\end{equation}
where $f\in\mL^{2}(\partial\Omega):= \{ v:\Omega\to \CC,\;\Vert
v\Vert_{\mL^{2}(\partial\Omega)}^2 := \int_{\partial\Omega}\vert v\vert^2 d\sigma
<+\infty\}$ is any square integrable function and
$\partial_{\bn}u:=\bn\cdot\nabla u$ with $\bn$ the vector field normal
to the boundary $\partial\Omega$ directed toward the exterior. The
wave number is modelled as a real constant $\kappa>0$. Following
widespread notations, we have considered the Sobolev space
$\mH^1(\Omega):= \{ v\in\mL^{2}(\Omega), \; \nabla
v\in\mL^{2}(\Omega)^d\}$ equipped with $\Vert
v\Vert_{\mH^1(\Omega)}^2:=\Vert \nabla v\Vert_{\mL^{2}(\Omega)}^2 +
\kappa^{2}\Vert v \Vert_{\mL^{2}(\Omega)}^2$.  Problem
\eqref{InitPb} can be put in the variational form: find $u\in
\mH^{1}(\Omega)$ such that $a(u,v) = \ell(v)\;\forall v\in
\mH^{1}(\Omega)$ where
\begin{equation}\label{VariationnalForm}
  \begin{array}{rl}
    a(u,v)  & \!\!\!:= \int_{\Omega}\nabla u\nabla \overline{v}
    - \kappa^{2}u\overline{v}\,d\bx - i\kappa\int_{\partial\Omega}u\overline{v}\,d\sigma\\[5pt]
    \ell(v) & \!\!\!:= \int_{\partial\Omega}f\overline{v}\,d\sigma.
  \end{array}
\end{equation}
Next we consider a regular triangulation $\mathcal{T}_h(\Omega)$ of
the computational domain $\overline{\Omega} = \cup_{\tau\in
  \mathcal{T}_h(\Omega)} \overline{\tau}$ and we denote
$\Vh(\Omega):=\{v\in\mathscr{C}^{0}(\overline{\Omega}):
v\vert_{\tau}\in \mathbb{P}_k(\overline{\tau})\;\forall \tau\in
\mathcal{T}_h(\Omega) \}\subset \mH^{1}(\Omega)$ a space of
$\mathbb{P}_k$-Lagrange finite element functions constructed on this
mesh, where $\mathbb{P}_k(\overline{\tau}):=\{$ polynomials of order
$\leq k$ on $\overline{\tau}$ $\}$.The associated discrete variational
formulation then writes
\begin{equation}\label{DiscreteVF}
  \begin{aligned}
    & \text{Find}\;u_h\in\Vh(\Omega)\;\text{such that}\\
    & a(u_h,v_h) = \ell(v_h)\;\; \forall v_h\in \Vh(\Omega).
  \end{aligned}
\end{equation}
Problem \eqref{DiscreteVF} shall be assumed to admit a unique
solution, which is simply equivalent to assuming that the
corresponding matrix (for a given choice of shape functions) is
invertible. The domain decomposition strategy that we
subsequently discuss aims at computing this solution.

\section{Decomposition of the computational domain}
In the perspective of domain decomposition, we need to introduce a
geometric decompositon of the computational domain. The strategy we
wish to consider belongs to the class of substructuring methods and
thus requires a non-overlapping partition of the computational domain
\begin{equation}
  \begin{aligned}
    \overline{\Omega}
    & = \overline{\Omega}_1\cup \dots\cup \overline{\Omega}_\mJ
    & \text{with} & \;\;\Omega_j\cap \Omega_k = \emptyset\;\text{for}\; j\neq k,\\
    \Sigma
    & = \Gamma_1\cup \dots\cup \Gamma_\mJ
    & \text{where} & \;\;\Gamma_j := \partial\Omega_j. 
  \end{aligned}
\end{equation}
Each $\Omega_j\subset \Omega$ will be a polyhedral assumed to be
exactly resolved by the triangulation i.e. $\overline{\Omega}_j =
\cup_{\tau\in \mathcal{T}_h(\Omega_j)}\overline{\tau}$ where
$\mathcal{T}_h(\Omega_j):=\{\tau\in \mathcal{T}_h(\Omega): \tau\subset
\Omega_j \}$. Following usual parlance, we shall call $\Sigma$ the
skeleton of the partition.

In practice the geometric decomposition above is obtained by means of
a graph partitionner. Such a decomposition a priori involves
cross-points i.e. points of adjacency of either three sub-domains or
two sub-domains meeting and the exterior boundary. The set of
cross-points is also refered to as wire basket in DDM related
litterature.  A major advantage of the DDM strategy we will consider
is its ability to handle cross-points properly.

\quad\\
We consider finite element spaces local to each subdomain
$\Vh(\Omega_j):=\{ v_h\vert_{\Omega_j}: v_h\in\Vh(\Omega)\}$, as well
as finite element spaces on local boundaries $\Vh(\Gamma_j):=\{
v_h\vert_{\Gamma_j}: v_h\in\Vh(\Omega)\}$. We shall also refer to
finite element functions defined on the skeleton
\begin{equation}
  \Vh(\Sigma):=\{ v_h\vert_{\Sigma}:\;v_h\in \Vh(\Omega)\}.
\end{equation}
We also need to consider volume based finite element functions that
are only piecewise continuous, with possible jumps through interfaces.
Such a space is naturally identified with a cartesian product. This
leads to setting
\begin{equation}\label{MultiVolumeSpaces}
\mbVh(\Omega):= \Vh(\Omega_1)\times \dots\times \Vh(\Omega_\mJ).
\end{equation}
We are interested in domain decomposition where behaviour of functions
at interfaces play a crucial role, so we also need to introduce
a space of traces at local boundaries $\mbVh(\Sigma)$
and the corresponding trace map $\Bh:\mbVh(\Omega)\to \mbVh(\Sigma)$
defined by
\begin{equation}\label{MultiTraceSpaces}
  \begin{aligned}
    & \mbVh(\Sigma):= \Vh(\Gamma_1)\times \dots\times \Vh(\Gamma_\mJ)\\
    & \Bh(v):= (v_1\vert_{\Gamma_{1}},\dots,v_\mJ\vert_{\Gamma_{\mJ}}).
  \end{aligned}
\end{equation}
Finally we need to embbed the space of trace on the skeleton
into the space of traces on local boundaries by means of a
restriction operator $\mR:\Vh(\Sigma)\to \mbVh(\Sigma)$ defined
subdomain-wise through
\begin{equation}\label{RestrictionOperator}
  \mR(v):=(v\vert_{\Gamma_1},\dots,v\vert_{\Gamma_\mJ}).
\end{equation}
The geometric decomposition that we have introduced above induces a
decomposition of the sesquilinear form \eqref{VariationnalForm} and leads
to an elementary reformulation of the discrete problem \eqref{DiscreteVF}. 
Define $\Ah:\mbVh(\Omega)\to \mbVh(\Omega)^*$ and $\boldsymbol{l}\in\mbVh(\Omega)^*$ by
\begin{equation}
  \begin{aligned}
    \langle \Ah(\bu),\bv\rangle
    & := \sum_{j=1,\dots,\mJ}\int_{\Omega_j}\nabla u_j\nabla v_j
    - \kappa^{2}u_jv_j\,d\bx - i\kappa
    \int_{\partial\Omega\cap\partial\Omega_j}u_jv_j d\sigma\\    
    \langle \bl,\bv\rangle
    & := \sum_{j=1,\dots,\mJ}\int_{\partial\Omega\cap\partial\Omega_j}f\,v_j\,d\bx 
  \end{aligned}
\end{equation}
for any $\bu = (u_1,\dots,u_{\mJ}), \bv =
(v_1,\dots,v_{\mJ})\in\mbVh(\Omega)$.  The operator $\Ah$ is block
diagonal, with each block associated with a different subdomain.
The initial discrete variational formulation can be
rewritten as follows: a function $u_h\in \Vh(\Omega)$ solves
\eqref{DiscreteVF} if and only if $\bu =
(u_h\vert_{\Omega_1},\dots,u_h\vert_{\Omega_\mJ})$ solves
\begin{equation}\label{SecondDiscreteFormulation}
  \begin{aligned}
    & \bu\in \mbXh(\Omega):=\{(v\vert_{\Omega_1},\dots,v\vert_{\Omega_{\mJ}}), v\in \Vh(\Omega)\}\\
    &  \text{and}\;\; \Ah(\bu),\bv\rangle = \langle \bl,\bv\rangle\;\;
    \forall\bv\in\mbXh(\Omega).
  \end{aligned}
\end{equation}

\section{Exchange operator}

To obtain a domain decomposition method, we need to further transform
\eqref{SecondDiscreteFormulation}. Our final formulation will be posed
in the function space $\mbVh(\Sigma)$ attached to the skeleton, and we
need to introduce a scalar product for this space i.e. an hermitian positive
definite operator
\begin{equation}\label{DefImpedance}
  \begin{aligned}
    & \mT:\mbVh(\Sigma)\to \mbVh(\Sigma)^*\;\; \text{such that}\\    
    & \mT = \mT^*\;\; \text{and}\;\; \langle \mT(\bv),\overline{\bv}\rangle>0\;
    \forall\bv\in\mbVh(\Sigma)\setminus\{0\}.
  \end{aligned}
\end{equation}
The operator $\mT^{-1}:\mbVh(\Sigma)^{*}\to \mbVh(\Sigma)$ induces a
scalar product on $\mbVh(\Sigma)^{*}$. This scalar product will be
used to quantify convergence of our domain decomposition algorithm. In
the subsequent analysis, the space $\mbVh(\Sigma)$
(resp. $\mbVh(\Sigma)^{*}$) will be equipped with the norm
\begin{equation}
  \begin{aligned}
    & \Vert \bv\Vert_{\mT}^2 &\hspace{-0.4cm}&:=\langle \mT(\bv),\overline{\bv}\rangle\\
    & \Vert \bq\Vert_{\mT^{-1}}^2&\hspace{-0.4cm}&:=\langle \mT^{-1}(\bq),\overline{\bq}\rangle\\
  \end{aligned}
\end{equation}
Based on the scalar products above, we are going to consider so-called
exchange operators that are responsible for enforcing transmission
conditions through interfaces and hence coupling between
subdomains. We define the exchange operator $\Pi:\mbVh(\Sigma)^{*}\to
\mbVh(\Sigma)^{*}$ by identity
\begin{equation}\label{DefExchg}
  \Pi:= 2\mT\mR(\mR^{*}\mT\mR)^{-1}\mR^* -\Id.
\end{equation}
Because $\mT\mR(\mR^{*}\mT\mR)^{-1}\mR^*$ is a $\mT^{-1}$-orthogonal
projector, it is clear that $\Pi$ is unitary by construction i.e.
$\Vert \Pi(\bp)\Vert_{\mT^{-1}} = \Vert \bp\Vert_{\mT^{-1}}\forall
\bp\in\mbVh(\Sigma)^*$. Besides it is a priori non-local in the sense
that it may couple trace function attached to distant subdomains.

\section{Skeleton formulation}
Now we describe a formulation equivalent to \eqref{SecondDiscreteFormulation} that serves
as the master equation of our domain decomposition algorithm. It will be posed on the
skeleton of the decomposition. In this skeleton equation, wave problems local to each subdomain
are written by means of a so-called scattering operator $\mS:\mbVh(\Sigma)^*\to\mbVh(\Sigma)^*$
defined by
\begin{equation}
  \begin{aligned}
    & \mS := \Id + 2i\mT\Bh(\Ah-i\Bh^*\mT\Bh)^{-1}\Bh^*\\
    & \bg := 2i\mT\Bh(\Ah-i\Bh^*\mT\Bh)^{-1}\bl
  \end{aligned}
\end{equation}
Since $\Ah,\Bh$ are both subdomain-wise block-diagonal, when $\mT$
is subdomain-wise block-diagonal, so is the scattering operator $\mS$.
This makes such an operator adapted to parallelism.
The following Proposition was established in \cite[Sec.6]{claeys2021nonself}.
\begin{prop}\quad\\
  The tuple function $\bu\in\mbXh(\Omega)$ solves
  \eqref{SecondDiscreteFormulation} if and only if there exists
  $\bq\in \mbVh(\Sigma)^*$ satisfying $\Bh^*\bq =
  (\Ah-i\Bh^*\mT\Bh)\bu - \bl$ and the skeleton formulation
  \begin{equation}\label{SkeletonFormulation}
    \begin{aligned}
      & \bq\in \mbVh(\Sigma)^*\;\text{and}\\
      & (\Id + \Pi\mS)\bq = \bg
    \end{aligned}
  \end{equation}
\end{prop}

\quad\\
From the proposition above, we see that if the skeleton formulation
\eqref{SkeletonFormulation} is solved, then the complete solution $\bu
= (\Ah-i\Bh^*\mT\Bh)^{-1}(\Bh^*\bq+\bl)$ can be reconstructed, and
this reconstruction step is fully parallel if the impedance operator
$\mT$ is block-diagonal.

\quad\\
A key feature of the skeleton formulation above is its strong coercivity
with respect to the scalar product induced by $\mT^{-1}$,
in spite of the a priori sign indefiniteness of \eqref{SecondDiscreteFormulation}.
The following result was established in \cite[Cor. 6.2]{claeys2021nonself}.
\begin{prop}\label{CoercivitySkeletonOperator}\quad\\
  The operator $\Id+\Pi\mS:\mbVh(\Sigma)^{*}\to \mbVh(\Sigma)^{*}$ is
  a bijection that fulfills a coercivity bound in the scalar product
  induced by $\mT^{-1}$. For all $\bq\in \mbVh(\Sigma)^*$ we have
  \begin{equation*}
    \begin{aligned}
      & \Re e\{\langle
      (\Id+\Pi\mS)\bq,\mT^{-1}\overline{\bq}\rangle\}\geq
      \frac{\gamma_{h}^{2}}{2}\Vert \bq\Vert_{\mT^{-1}}^{2}\\
      & \text{where}\;\; \gamma_h:=\mathop{\inf}_{\bq\in
        \mbVh(\Sigma)^*\setminus\{0\}}
      \Vert (\Id +\Pi\mS)\bq\Vert_{\mT^{-1}}/\Vert \bq\Vert_{\mT^{-1}}
    \end{aligned}
  \end{equation*}
\end{prop}
\noindent 
This result garantees a good behaviour of classical iterative solvers
such as Richardson or GMRes when applied to the skeleton formulation
\eqref{SkeletonFormulation}.  The coercivity estimate above leads to
convergence bounds for such solvers, see e.g.
\cite{zbMATH05029264}.

\quad\\
In a Krylov solver, the matrix-vector operation $\bv\mapsto
(\Id+\Pi\mS)\bv$ is a crucial step that has enormous impact on the
computational cost of the overall solution procedure. In a situation
where the impedance operator $\mT$ is block-diagonal, the
matrix-vector product $\bv\mapsto \mS\bv$ is ambarassingly parallel as
the scattering operator is itself block-diagonal.

As a consequence, the cost of the operation $\bv\mapsto \Pi(\bv)$ is
critical.  This reduces to performing $\bv\mapsto
(\mR^*\mT\mR)^{-1}\bv$ which, if performed in a naive way (like with a
direct Cholesky solver), might be costly. Indeed solving a linear
system attached to $\mR^*\mT\mR$ is required each time a matrix-vector
product is required inside the global Krylov solver.

\quad\\
The main goal of the present article is to show how this core step,
that is at the center of our DDM strategy, can be optimized. We wish
to exhibit how each such linear solve can benefit from the previous 
solves by means of a simple recycling strategies in such a way that the
convergence speed of the overall DDM algorithm is not deteriorated.

\quad\\
For the sake of concreteness, we consider a Richardson
solver, see e.g. Example 4.1 in \cite{zbMATH01953444}, as a model
iterative solver for the skeleton formulation
\eqref{SkeletonFormulation}.  Starting from a trivial initial guess
$\bq^{(0)} = 0$ and choosing $\alpha\in (0,1)$ as relaxation
parameter, Richardson's iteration takes the form
\begin{equation}\label{RichardsonInit}
  \bq^{(n+1)} = (1-\alpha)\bq^{(n)}-\alpha\Pi\mS\bq^{(n)}+\alpha \bg.
\end{equation}
Taking account of the expression of the exchange operator given by
\eqref{DefExchg}, this iteration can be decomposed as follows
\begin{equation}\label{RichardsonIteration}
  \begin{array}{l}
    \textbf{Exact Richardson iteration}\\[5pt]
    \begin{aligned}    
      \bp^{(n+1)} & = (\mR^*\mT\mR)^{-1}\mR^*\mS\bq^{(n)}\\
      \bq^{(n+1)} & = ((1-\alpha)\Id + \alpha\mS)\bq^{(n)}-2\alpha\mT\mR(\bp^{(n+1)})+\alpha \bg
    \end{aligned}
  \end{array}
\end{equation}
We dub this algorithm "exact" because all steps in this procedure
are assumed to be conducted without any error. In particular, it is
assumed that no approximation is made when evaluating the action of
$(\mR^*\mT\mR)^{-1}$ in the first line above.

\section{Approximation of the exchange operator}
In a large distributed memory environment, the linear solve
of $\bq\mapsto(\mR^*\mT\mR)^{-1}\bq$ cannot be achieved exactly but
should rely instead on an approximate PCG solve that consists in
Galerkin projections onto a Krylov space that depends itself on the
right hand side, see e.g. \cite[\S 2.3]{MR3024841}.

\quad\\
Consider a linear map $\Prec:\Vh(\Sigma)^*\to \Vh(\Sigma)$ that, in the
subsequent analysis, shall play the role of a preconditioner. For a
given integer $k\geq 1$ and an initial guess $\bx_0\in\Vh(\Sigma)$,
define the $k$-th order Krylov space
\begin{equation}
  \begin{aligned}
    & \mathscr{K}_k(\bb,\bx_0):=
    \textrm{vect}\{\br_0,(\Prec\mR^*\mT\mR)\br_0,\dots,(\Prec\mR^*\mT\mR)^{k-1}\br_0\}\\
    & \text{where}\;\;\br_0:=\Prec\bb-(\Prec\mR^*\mT\mR)\bx_0.
  \end{aligned}
\end{equation}
Define $\mrm{PCG}_k:\Vh(\Sigma)\times\Vh(\Sigma)^*\to
\Vh(\Sigma)$ as the map that takes a right-hand side
$\bb\in\Vh(\Sigma)^*$ and an initial guess $\bx_0\in\Vh(\Sigma)$, and
returns the result of $k$ steps of a PCG algorithm preconditionned
with $\Prec$, see \cite[\S 9.2]{zbMATH01953444}. Following the
interpretation of Krylov methods in terms of projections, see
\cite[Chap.5\,\&\,6]{zbMATH01953444} and \cite[Chap.2]{MR3024841}, it
is caracterised as the unique solution to the following minimization problem
\begin{equation}\label{CaracPCG}
  \begin{aligned}
    & \mrm{PCG}_k(\bx_0,\bb)\in \bx_0+\mathscr{K}_k(\bx_0,\bb)\quad \text{and}\\ 
    & \Vert (\mR^*\mT\mR)^{-1}\bb - \mrm{PCG}_k(\bx_0,\bb) \Vert_{\mR^*\mT\mR} =
    \min_{\eta\in \bx_0+\mathscr{K}_k(\bx_0,\bb)}\Vert (\mR^*\mT\mR)^{-1}\bb - \eta \Vert_{\mR^*\mT\mR}.
  \end{aligned}
\end{equation}
Let us see how the preconditioned conjugate gradient map $\mrm{PCG}_k$
can be inserted into Richardson's iteration
\eqref{RichardsonIteration}. A naive approach would systematically
take $\bx_0 = 0$ as initial guess which leads to considering the
relation $\bp^{(n+1)} = \mrm{PCG}_k(0,\mR^*\mS\bq^{(n)})$ for the
first equation in \eqref{RichardsonIteration}. However, because the
overall DDM algorithm is supposed to converge, the right-hand sides
$\mR^*\mS\bq^{(n)}$ should themselve form a converging sequence. As a
consequence $\mR^*\mS\bq^{(n-1)}$ should remain close to
$\mR^*\mS\bq^{(n)}$, so that taking $\bx_0 = \bp^{(n)}$ appears as a
natural recycling strategy. The modified DDM strategy then takes the
form
\begin{equation}\label{InexactRichardson}
  \begin{array}{l}
    \textbf{Approximate Richardson iteration}\\[5pt]
    \begin{aligned}
      \tilde{\bp}^{(n+1)} & = \mrm{PCG}_k(\tilde{\bp}^{(n)},\mR^*\mS\tilde{\bq}^{(n)})\\
      \tilde{\bq}^{(n+1)} & = ((1-\alpha)\Id + \alpha\mS)\tilde{\bq}^{(n)}
      -2\alpha\mT\mR(\tilde{\bp}^{(n+1)})+\alpha \bg
    \end{aligned}
  \end{array}
\end{equation}
The parameter $k$ that represents the dimension of the Krylov space
is assumed fixed and independent of $n$. Because the map $\mrm{PCG}_k$ is
not linear, we underline that the iterative procedure above is itself
non-linear.

\section{Numerical experiment}

We now present a numerical experiment illustrating the strategy
described above. We will consider Problem \eqref{InitPb} set in $\RR^d
= \RR^2$ in the computational domain $\Omega = (-1,+1)^{2}\setminus
(-0.25,+0.25)^{2}$ represented below in Fig.\ref{Fig1}. We take
$\lambda = 1/10$ hence a wavenumber $\kappa = 2\pi/\lambda \simeq
62.83$, and the source function $f(\bx) = \exp(i\kappa\,\bd\cdot\bx)$
with $\bd = (-1/\sqrt{2}, +1/\sqrt{2},0)$.  The problem is discretized
with $\mathbb{P}_1$-Lagrange finite elements based on a mesh generated
with \textsc{gmsh} and partitioned with \textsc{metis} in 16
subdomains. The computations were run sequentially with the C++
library \textsc{ddmtool} on a laptop with
Intel\textsuperscript{\textregistered} Core\texttrademark \,i7-1185G7
processor with 62.5 Gb of RAM.

\quad\\
All the numerical experiments of the present section have been
conducted with the same fixed mesh that contains 175794 nodes (which
is also the dimension of $\Vh(\Omega)$) and 349588 triangles. The
maximum mesh element size equals $h = 0.005$ so this discretization
represents approximately $\lambda/h = 20$ points per wavelength. We
compute a reference solution of \eqref{SkeletonFormulation} by means
of a direct solver using \textsc{umpfpack}.  This reference solution
will be denoted $\bq^{(\infty)}$ subsequently. Its real part is
plotted in the right hand side of Figure \ref{Fig1} below.
\begin{figure}[H]
  \begin{center}
    \includegraphics*[height=5cm,trim= 0cm 2cm 0cm 1cm]{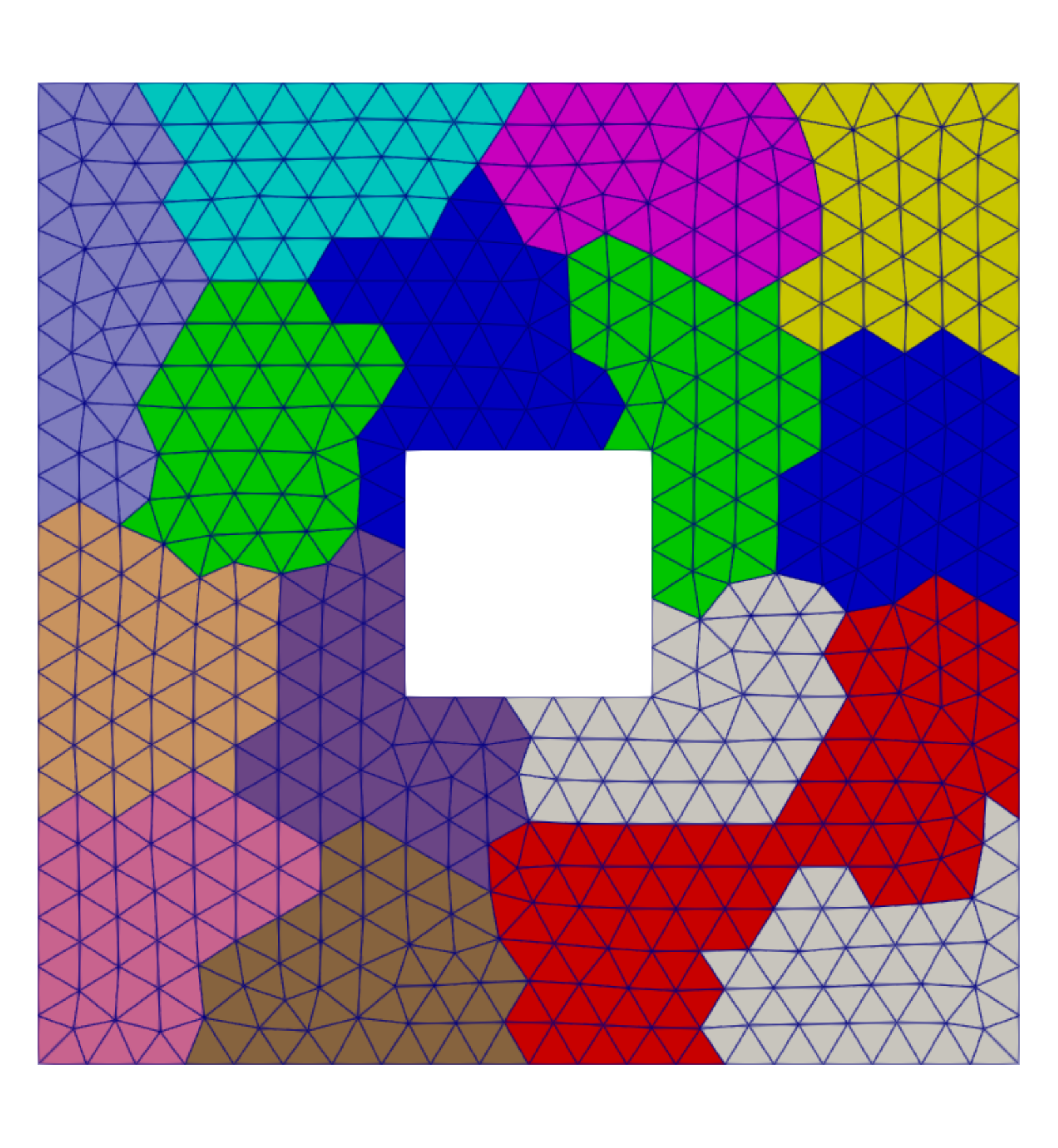}
    \hspace{1cm}
    \includegraphics*[height=5cm,trim= 0cm 2cm 0cm 1cm]{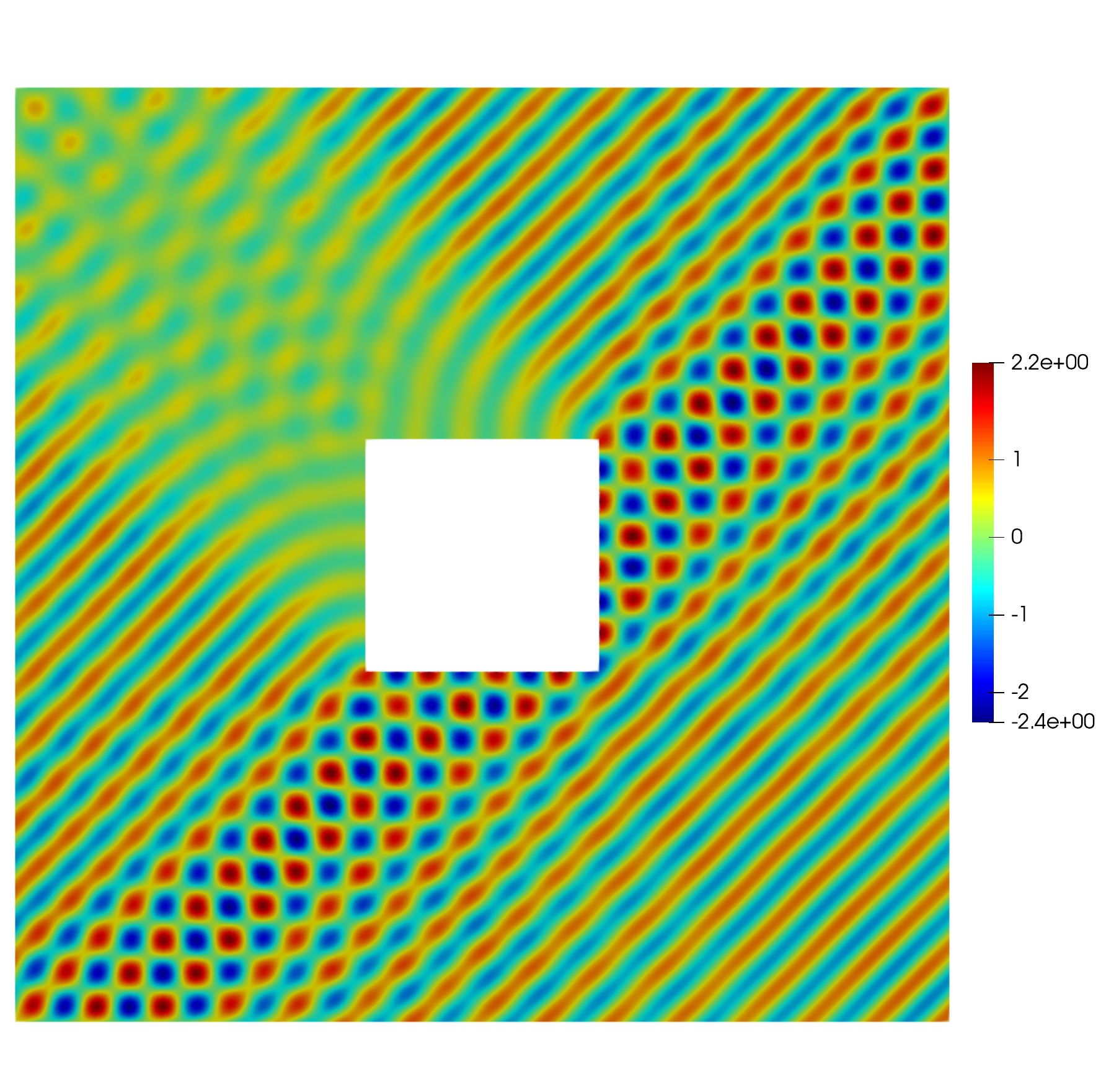}
  \end{center}
  \caption{Left: Computational domain.
    Right: Real part of the reference solution.}\label{Fig1}
\end{figure}
\noindent
We construct the impedance operator $\mT =
\mrm{diag}(\mT_1,\dots,\mT_\mJ):\mbVh(\Sigma)\to \mbVh(\Sigma)^*$
following the strategy advocated in
\cite[Chap.8]{parolin:tel-03118712}, \cite{ParolinCollinoJoly},
\cite[\S4.2]{MR4507159}. In the present case, this boils down to
selecting a subset of each subdomain $\tilde{\Omega}_j\subset
\Omega_j$ consisting in 5 layers of elements neighbouring $\Gamma_j =
\partial\Omega_j$ (in particular $\Gamma_j\subset
\partial\tilde{\Omega}_j$) and defining each
$\mT_j:\Vh(\Gamma_j)\to\Vh(\Gamma_j)^{*}$ as the unique hermitian
positive definite linear map satisfying the minimization property
\begin{equation*}
    \langle \mT_j(v),v\rangle := \min\{\; 
    \Vert \tilde{v}\Vert_{\mL^{2}(\tilde{\Omega}_j)}^2+\kappa^2
    \Vert \tilde{v}\Vert_{\mL^{2}(\tilde{\Omega}_j)}^2+ \kappa
    \Vert \tilde{v}\Vert_{\mL^{2}(\partial\tilde{\Omega}_j\setminus\Gamma_j)}^2,
    \;\tilde{v}\in \Vh(\tilde{\Omega}_j),\;\tilde{v}\vert_{\Gamma_j} = v\}.
\end{equation*}
The actual evaluation of this impedance operator rests on a (sparse)
Cholesky factorization performed by means of \textsc{umfpack} locally
in each $\tilde{\Omega}_j$.  As for the preconditioner $\Prec$ for the
linear solve associated to the operation $\bb\mapsto
(\mR^*\mT\mR)^{-1}\bb$, we choose the single level Neumann-Neumann
preconditioner, see \cite[\S 7.8.1]{MR3450068} or \cite[\S
  6.2]{MR2104179}.

\quad\\
In a first experiment we run a variant of Algorithm
\eqref{InexactRichardson} with $\alpha = 1/2$, where the initial guess
of PCG is chosen trivial i.e. we set $\tilde{\bp}^{(n+1)} =
\mrm{PCG}_k(0, \bb)$ where $\bb = \mR^*\mS\,\tilde{\bq}^{(n)}$, and
the PCG algorithm is executed until the relative residual
error $\Vert\bb-(\mR^*\mT\mR)\mrm{PCG}_k(0,\bb)\Vert_{\Prec}/
\Vert\bb\Vert_{\Prec}\leq 1e-20$ is reached. In Figure \ref{Fig2} below, we plot
the norm of the error $\Vert \bq^{(\infty)} -
\tilde{\bq}^{(n)}\Vert_{\mT^{-1}}/\Vert \bq^{(\infty)}\Vert_{\mT^{-1}}$ of
Algorithm \eqref{InexactRichardson} versus the iteration number $n$.
\begin{figure}[H]
  \begin{center}
    \includegraphics*[height=5cm,trim= 0cm 0.25cm 0cm 0.5cm]{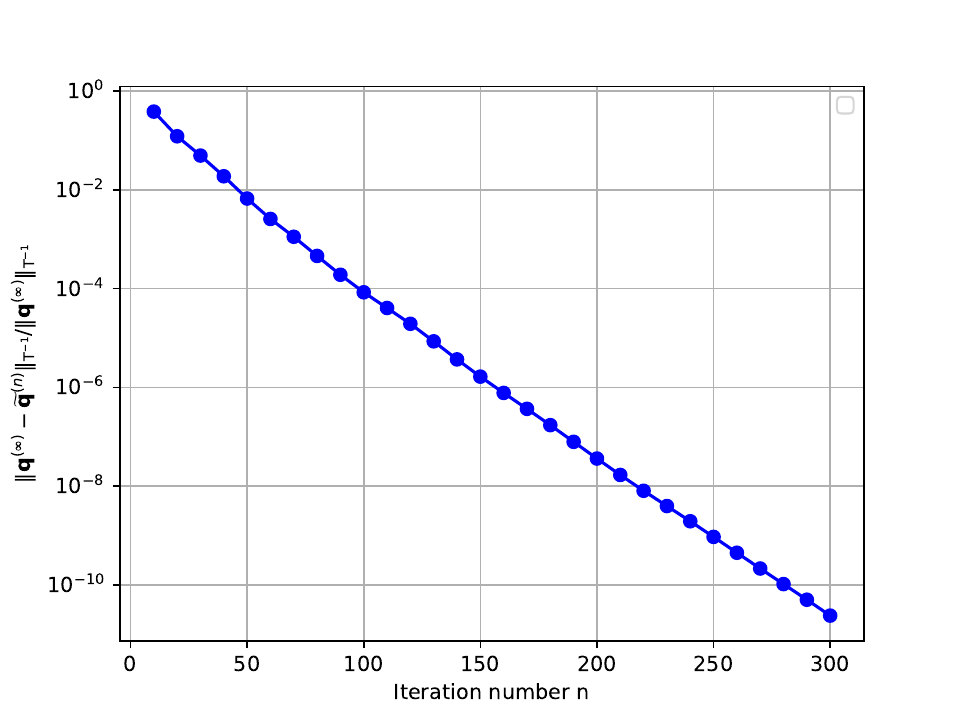}
    \includegraphics*[height=5cm,trim= 0cm 0.25cm 0cm 0.5cm]{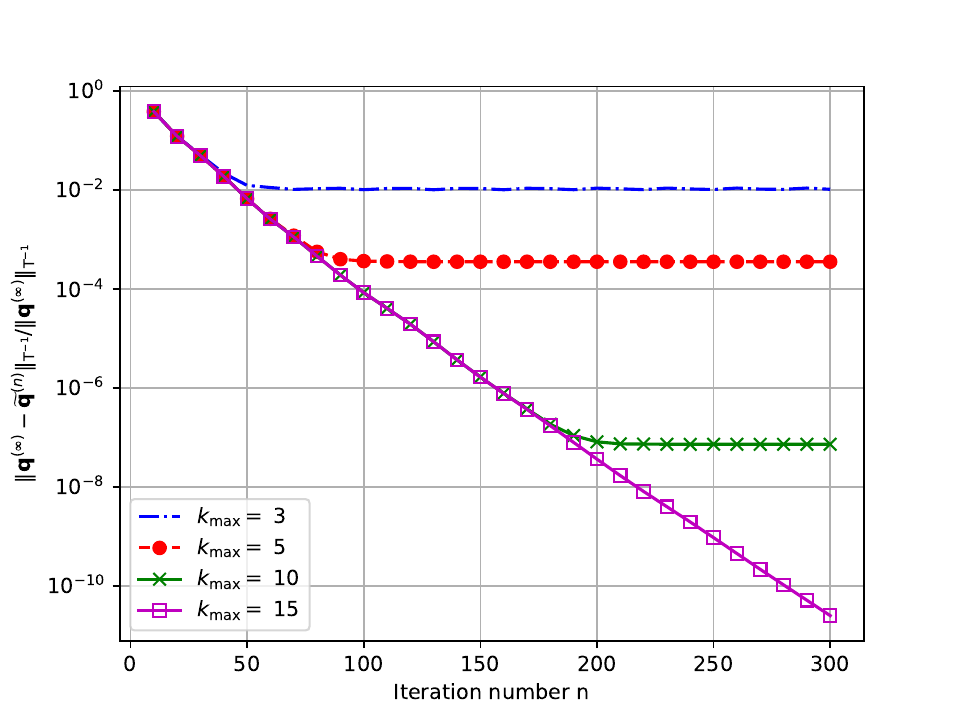}
  \end{center}
  \caption{Relative error $\Vert \bq^{(\infty)} -
    \tilde{\bq}^{(n)}\Vert_{\mT^{-1}}/\Vert
    \bq^{(\infty)}\Vert_{\mT^{-1}}$ versus iteration number $n$ in
    Richardson's Algorithm. No initial guess of PCG i.e. no recycling strategy.
    Left: no limitation on PCG iteration count i.e. $k_{\textsc{max}} =
    \infty$. Right: imposing in addition $k\leq
    k_{\textsc{max}}$.}\label{Fig2}
\end{figure}
\noindent 
On the left picture of Fig.\ref{Fig2}, we take as many iterations of
PCG as needed. For this plot, as a matter of fact, 14 iterations of PCG
take place at each iteration of the global Richardson algorithm. On
the right hand side of Fig.\ref{Fig2}, we plot the same graph except
that this time the number of iterations of PCG is limited $k\leq
k_{\textsc{max}}$ for several values of $k_{\textsc{max}}$. We see
that when the number of PCG iterations is limited, the error $\Vert
\bq^{(\infty)} - \tilde{\bq}^{(n)}\Vert_{\mT^{-1}}/\Vert
\bq^{(\infty)}\Vert_{\mT^{-1}}$ of the global Richardson algorithm
decays normally until it reaches a certain critical value where it
stalls. This plateau phenomenon appears here for $k_{\textsc{max}}< 14$.

\quad\\
Next we launch the same computation, again limiting the number of PCG
iterations $k\leq k_{\textsc{max}}$. This time though, we choose the
initial guess from the previous iterate as described before
i.e. $\tilde{\bp}^{(n+1)} =\mrm{PCG}_k(\bx_0, \bb)$ with $\bb =
\mR^*\mS\,\tilde{\bq}^{(n)}$ and $\bx_0 = \tilde{\bp}^{(n)}$.  We do
this for the two values $k_{\textsc{max}} = 5$ and $k_{\textsc{max}} =
10$ and plot the corresponding error $\Vert \bq^{(\infty)}
-\tilde{\bq}^{(n)}\Vert_{\mT^{-1}}/\Vert
\bq^{(\infty)}\Vert_{\mT^{-1}}$ versus $n$. This time the error keeps
on decaying without reaching any plateau.
\begin{figure}[H]
  \begin{center}
    \includegraphics*[height=5cm,trim= 0.cm 0.25cm 0cm 0.cm]{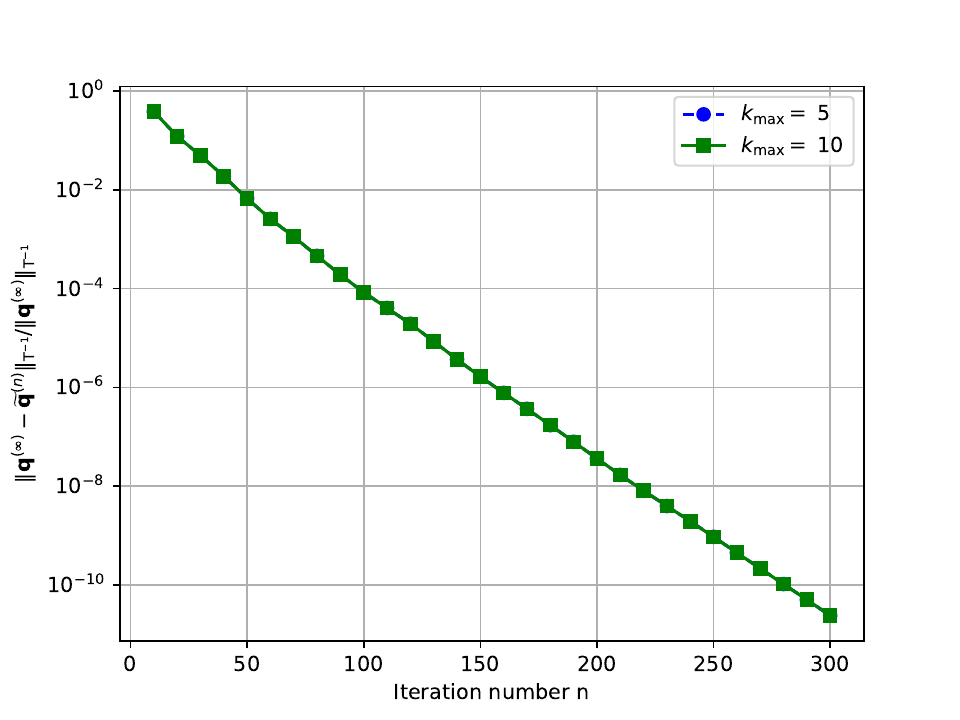}
  \end{center}
  \caption{Relative error $\Vert \bq^{(\infty)} -
    \tilde{\bq}^{(n)}\Vert_{\mT^{-1}}/\Vert
    \bq^{(\infty)}\Vert_{\mT^{-1}}$ versus iteration number $n$ in
    Richardson's Algorithm with a recycled initial guess for PCG
    and $k\leq k_{\textsc{max}} = 5,10$.}\label{Fig3}
\end{figure}
\noindent 
The plot of Fig.\ref{Fig3} looks identical to the one in the left hand
side of Fig.\ref{Fig2}. This suggests that, when combining a
truncation of PCG with the simple recycling strategy described above,
the error decay of the global Richardson algorithm does not experience
any deterioration.  Keeping this strategy consisting in both recycling
and truncating PCG, in the next table, we give the number of
iterations $n$ required to reach a relative tolerance $\Vert
\bq^{(\infty)} -\tilde{\bq}^{(n)}\Vert_{\mT^{-1}}/\Vert\bq^{(\infty)}\Vert_{\mT^{-1}}<1e-10$.
\begin{center}
  \begin{tabular}{l|ccccccc}
    $k_{\textsc{max}}$
    & 20  & 15  & 10   & 5   & 3   & 2   & 1\\\hline
    $\#$iter
    & 281 & 281 & 281  & 281 & 282 & 286 & 616
  \end{tabular}
\end{center}
This clearly indicates that, when using a recycling strategy, only a
few PCG iterations are needed to maintain the same convergence speed
for the overall Richardson algorithm. This represents a clear
computational gain since the cost of each iteration of
\eqref{InexactRichardson} depends directly on $k_{\textsc{max}}$.

\quad\\
The previous result shows that, although the operator $\Pi$ is
non-local, in practice its action may be evaluated with a cost
corresponding to a small number of matrix-vector products from the
operator $\mT$. In the next section, we provide theoretical analysis
supporting this conclusion.

\section{Convergence analysis}
We exhibited an efficient heuristic to approximate the action of the
non-local exchange operator i.e. the first line in
\eqref{RichardsonIteration}. It consists in combining a brutal
truncation of PCG with a basic recycling scheme. We provided numerical
evidence supporting the relevance of this strategy.

We seek now to obtain a theoretical explanation for the performance of
this approach. Instead of trying to systematically derive the sharpest
estimates, at certain points of our analysis we will take upper bounds
that are larger than strictly required which, hopefully, will help
simplify the calculus. To begin with, we re-arrange the approximate Richardson
iteration \eqref{InexactRichardson}, 
\begin{equation}\label{FirstRecursion}
  \begin{aligned}
    & \tilde{\bp}^{(n+1)} -\tilde{\bp}^{(n+1)}_{\infty}
    = \mrm{PCG}_k(\tilde{\bp}^{(n)},\mR^*\mS\tilde{\bq}^{(n)}) -
    (\mR^*\mT\mR)^{-1}\mR^*\mS\tilde{\bq}^{(n)}\\
    & \tilde{\bq}^{(n+1)} = ((1-\alpha)\Id - \alpha\Pi\mS)\tilde{\bq}^{(n)}
    -2\alpha\mT\mR(\tilde{\bp}^{(n+1)}-\tilde{\bp}^{(n+1)}_{\infty})+\alpha \bg\\[10pt]
    & \text{where}\quad   \tilde{\bp}_{\infty}^{(n)}:=
    (\mR^*\mT\mR)^{-1}\mR^*\mS\tilde{\bq}^{(n-1)}.
  \end{aligned}
\end{equation}
Focusing on the first line of \eqref{FirstRecursion}, we try to
estimate the decay of the left hand side, making use of the classical
convergence estimate for the conjugate gradient, see e.g. Corollary
5.6.7 in \cite{MR3024841} that, in our notations, yields the following
inequality
\begin{equation}\label{EstimateCG}
  \begin{aligned}
    & \Vert (\mR^*\mT\mR)^{-1}\bb - \mrm{PCG}_k(\bx_0,\bb)\Vert_{\mR^*\mT\mR} \leq
    \epsilon \Vert (\mR^*\mT\mR)^{-1}\bb - \bx_0\Vert_{\mR^*\mT\mR}\\[5pt]
    & \text{with}\;\; \epsilon =
    2\Big(\frac{
      \sqrt{\mrm{cond}(\Prec\mR^*\mT\mR)}-1}{
    (\sqrt{\mrm{cond}(\Prec\mR^*\mT\mR)}+1}\Big)^k
  \end{aligned}
\end{equation}
and $\mrm{cond}(\mL)$ refers to the spectral condition number i.e.
$\mrm{cond}(\mL) = \sup\mathfrak{S}(\mL)/\inf\mathfrak{S}(\mL)$ where
$\mathfrak{S}(\mL)$ is the spectrum of a linear map
$\mL:\Vh(\Sigma)\to \Vh(\Sigma)$. For the moment, we assume that
$k$ is chosen sufficiently large hence $\epsilon$ as small as required.
We shall come back an discuss later on the choice of the parameter $k$.
Injecting Estimate \eqref{EstimateCG} into \eqref{FirstRecursion} leads
to following inequality
\begin{equation*}
  \begin{aligned}
    \Vert \tilde{\bp}^{(n+1)}_{\infty} - \tilde{\bp}^{(n+1)}\Vert_{\mR^*\mT\mR}
    & \leq \epsilon \Vert \tilde{\bp}^{(n+1)}_{\infty} - \tilde{\bp}^{(n)}\Vert_{\mR^*\mT\mR}\\
    & \leq \epsilon \Vert \tilde{\bp}^{(n)}_{\infty} - \tilde{\bp}^{(n)}\Vert_{\mR^*\mT\mR}
    + \epsilon \Vert \tilde{\bp}^{(n+1)}_{\infty} - \tilde{\bp}^{(n)}_{\infty}\Vert_{\mR^*\mT\mR}.
  \end{aligned}
\end{equation*}
By the very definition of the auxiliary variable in \eqref{FirstRecursion}, we have
$\tilde{\bp}^{(n+1)}_{\infty} - \tilde{\bp}^{(n)}_{\infty} =
(\mR^*\mT\mR)^{-1}\mR^{*}\mS(\bw)$ where $\bw = \tilde{\bq}^{(n)} - \tilde{\bq}^{(n-1)}$.
On the other hand, since $\mT\mR(\mR^*\mT\mR)^{-1}\mR^*$ is a $\mT^{-1}$-orthogonal
projection, and $\mS$ is a contraction with respect to $\Vert\;\Vert_{\mT^{-1}}$
according to Lemma 5.2 in \cite{claeys2021nonself}, we have
$\Vert \tilde{\bp}^{(n+1)}_{\infty} - \tilde{\bp}^{(n)}_{\infty}\Vert_{\mR^*\mT\mR}
=\Vert \mT\mR(\mR^*\mT\mR)^{-1}\mR^{*}\mS(\bw)\Vert_{\mT^{-1}}
\leq \Vert \bw\Vert_{\mT^{-1}} = \Vert \tilde{\bq}^{(n)} - \tilde{\bq}^{(n-1)}\Vert_{\mT^{-1}}$.
From this we obtain
\begin{equation}\label{Estimate1}
  \begin{aligned}
    \Vert \tilde{\bp}^{(n+1)}_{\infty} - \tilde{\bp}^{(n+1)}\Vert_{\mR^*\mT\mR}
    \leq \epsilon \Vert \tilde{\bp}^{(n)}_{\infty} - \tilde{\bp}^{(n)}\Vert_{\mR^*\mT\mR}
    + \epsilon \Vert \tilde{\bq}^{(n)} - \tilde{\bq}^{(n-1)}\Vert_{\mT^{-1}}.
  \end{aligned}
\end{equation}
Coming back to the approximate Richardson iteration
\eqref{FirstRecursion}, we now focus on the second line. We take the
difference of two successive iterates, which yields
\begin{equation}
  \begin{aligned}
    \tilde{\bq}^{(n+1)}-\tilde{\bq}^{(n)}
    & = ((1-\alpha)\Id - \alpha\Pi\mS)(\tilde{\bq}^{(n)}-\tilde{\bq}^{(n-1)})\\
    & \phantom{=}
    -2\alpha\mT\mR(\tilde{\bp}^{(n+1)}-\tilde{\bp}^{(n+1)}_{\infty})
    +2\alpha\mT\mR(\tilde{\bp}^{(n)}-\tilde{\bp}^{(n)}_{\infty})
  \end{aligned}
\end{equation}
Next we bound the norm of the left-hand side above,
taking account of Inequality \eqref{Estimate1}, and introducing
the continuity modulus of $(1-\alpha)\Id - \alpha\Pi\mS$ with respect to the
norm induced by $\mT^{-1}$. This yields 
\begin{equation}\label{Estimate2}
  \begin{aligned}
    & \Vert \tilde{\bq}^{(n+1)}-\tilde{\bq}^{(n)}\Vert_{\mT^{-1}}\\
    & \leq \rho \Vert \tilde{\bq}^{(n)}-\tilde{\bq}^{(n-1)}\Vert_{\mT^{-1}}
    + 2\alpha  \Vert \tilde{\bp}^{(n+1)}_{\infty} - \tilde{\bp}^{(n+1)}\Vert_{\mR^*\mT\mR}
    + 2\alpha \Vert \tilde{\bp}^{(n)}_{\infty} - \tilde{\bp}^{(n)}\Vert_{\mR^*\mT\mR}\\
    & \leq (\rho+2\alpha\epsilon)\Vert \tilde{\bq}^{(n)}-\tilde{\bq}^{(n-1)}\Vert_{\mT^{-1}}
    + 2\alpha (1+\epsilon)\Vert \tilde{\bp}^{(n)}_{\infty} - \tilde{\bp}^{(n)}\Vert_{\mR^*\mT\mR}\\[10pt]
    & \text{where}\quad   \rho:=\sup_{\bq\in \mbVh(\Sigma)^*\setminus\{0\}}
    \frac{\Vert ((1-\alpha)\Id - \alpha\Pi\mS)\bq\Vert_{\mT^{-1}}}{\Vert \bq\Vert_{\mT^{-1}}}.
  \end{aligned}
\end{equation}
Now we can gather \eqref{Estimate1} and \eqref{Estimate2}
to form a system of $2\times 2$ recursive inequalities.
We slightly rescale the inequalities and apply the majorizations
$\epsilon \leq 2\epsilon$ and $\alpha\leq 1$ to make the analysis
more comfortable. Set
\begin{equation}
  \mathfrak{R} :=
  \left\lbrack\begin{array}{cc}
  \rho+2\epsilon & 2\sqrt{\epsilon(1+\epsilon)}\\
  2\sqrt{\epsilon(1+\epsilon)} & 2\epsilon
  \end{array}\right\rbrack,\quad
  \mathfrak{e}^{(n)}:=
  \left\lbrack\begin{array}{l}
  \Vert \tilde{\bq}^{(n)}-\tilde{\bq}^{(n-1)}\Vert_{\mT^{-1}}\\
  \Vert \tilde{\bp}^{(n)}_{\infty} - \tilde{\bp}^{(n)}\Vert_{\mR^*\mT\mR}\sqrt{1+1/\epsilon}
  \end{array}\right\rbrack.
\end{equation}
Given two vectors of positive coefficients $u,v\in \mathbb{R}_+^{2}$
with $u = (u_1,u_2)$ and $v = (v_1,v_2)$, we shall write $u\leq v\iff
u_j\leq v_j $ for $j=1,2$. With this notation we obtain
$\mathfrak{e}^{(n+1)}\leq
\mathfrak{R}\cdot\mathfrak{e}^{(n)}$. Iterating over $n$, since all
coefficients are positive, we obtain $\mathfrak{e}^{(n)}\leq
\mathfrak{R}^{n}\cdot\mathfrak{e}^{(0)}$ with $\mathfrak{e}^{(0)} =
(\Vert \tilde{\bq}^{(0)}\Vert_{\mT^{-1}},0)^{\top}$.  For $\bx =
(x_1,x_2)\in\RR^2$, denoting $\vert \bx\vert := (\vert
x_1\vert^2+\vert x_2\vert^2)^{1/2}$ and the associated matrix norm
$\vert \mathfrak{R}\vert = \sup_{\bx\in \mathbb{R}^2\setminus\{0\}}
\vert \mathfrak{R}\bx\vert/\vert\bx\vert$.

\begin{lem}\quad\\
  Under the condition that $\vert\mathfrak{R}\vert<1$, the sequence
  $\tilde{\bq}^{(n)}$ defined by \eqref{InexactRichardson} converges
  toward $\bq^{(\infty)} := (\Id+\Pi\mS)^{-1}\bg$ the solution to
  \eqref{SkeletonFormulation} and we have $\Vert
  \tilde{\bq}^{(n)}-\bq^{(\infty)}\Vert_{\mT^{-1}}\leq \Vert
  \tilde{\bq}^{(0)}\Vert_{\mT^{-1}}\, \vert
  \mathfrak{R}\vert^{n+1}/(1-\vert \mathfrak{R}\vert)$ for all $n\geq
  0$.
\end{lem}
\noindent \textbf{Proof:}

First of all observe that
$\Vert \tilde{\bq}^{(n)}-\tilde{\bq}^{(n-1)}\Vert_{\mT^{-1}}
\leq \vert \mathfrak{e}^{(n)}\vert\leq 
\vert\mathfrak{R}\vert^{n}\cdot\vert\mathfrak{e}^{(0)}\vert$.
and $\vert\mathfrak{e}^{(0)}\vert = \Vert \tilde{\bq}^{(0)}\Vert_{\mT^{-1}}$. 
Now pick arbitrary integers $n,m$ with $m>n$. Under the assumption that
$\vert \mathfrak{R}\vert<1$, we have the following estimate 
\begin{equation}\label{CauchyCriterion}
  \begin{aligned}
    \Vert\tilde{\bq}^{(m)} - \tilde{\bq}^{(n)}\Vert_{\mT^{-1}}
    & \leq \sum_{\nu = n+1}^{m}  \Vert\tilde{\bq}^{(\nu)} - \tilde{\bq}^{(\nu-1)}\Vert_{\mT^{-1}}
    \leq \Vert\tilde{\bq}^{(0)}\Vert_{\mT^{-1}} \sum_{\nu = n+1}^{m}\vert \mathfrak{R}\vert^{\nu}\\
    & \leq \Vert\tilde{\bq}^{(0)}\Vert_{\mT^{-1}}\vert \mathfrak{R}\vert^{n+1}/(1-\vert \mathfrak{R}\vert).
  \end{aligned}
\end{equation}
This proves that the sequence $\tilde{\bq}^{(n)}$ is of Cauchy type in
the norm $\Vert \cdot\Vert_{\mT^{-1}}$ and admits a limit that we
denote $\tilde{\bq}^{(\infty)}$. Letting $m\to \infty$ in
\eqref{CauchyCriterion} yields $\Vert\tilde{\bq}^{(n)} -
\tilde{\bq}^{(\infty)}\Vert_{\mT^{-1}}\leq
\Vert\tilde{\bq}^{(0)}\Vert_{\mT^{-1}}\vert
\mathfrak{R}\vert^{n+1}/(1-\vert \mathfrak{R}\vert)$, so there only
remains to prove that $\tilde{\bq}^{(\infty)} = \bq^{(\infty)}$.
Observe now that we also have
\begin{equation*}
\Vert \tilde{\bp}_{\infty}^{(n)}-
\tilde{\bp}^{(n)} \Vert_{\mR^*\mT\mR}
\leq \frac{\vert \mathfrak{e}^{(n)}\vert}{\sqrt{1+1/\epsilon}}
\leq \frac{\vert \mathfrak{e}^{(0)}\vert\,\vert\mathfrak{R}\vert^{n}}{
  \sqrt{1+1/\epsilon}}
\end{equation*}
which proves that $\Vert \tilde{\bp}_{\infty}^{(n)}- \tilde{\bp}^{(n)}
\Vert_{\mR^*\mT\mR} \to 0$.  Next, coming back to
\eqref{InexactRichardson}, recall that we have the relation
$\tilde{\bq}^{(n+1)} = ((1-\alpha)\Id - \alpha\Pi\mS)\tilde{\bq}^{(n)}
-2\alpha\mT\mR(\tilde{\bp}^{(n+1)}-\tilde{\bp}^{(n+1)}_{\infty})+\alpha
\bg$.  Taking $n\to \infty$ in the previous relation, we obtain that
$\tilde{\bq}^{(\infty)}$ satisfies the equation
$\tilde{\bq}^{(\infty)} = ((1-\alpha)\Id -
\alpha\Pi\mS)\tilde{\bq}^{(\infty)} + \alpha \bg$. After 
re-arrangement, this leads to $\tilde{\bq}^{(\infty)} =
(\Id+\Pi\mS)^{-1}\bg = \bq^{(\infty)}$, which ends the
proof.\hfill $\Box$

\quad\\
A remarkable conclusion that can be drawn from the previous lemma is
that, with recycling, truncation of the PCG algorithm for the
computation of the exchange operator $\Pi$ does not induce any
consistency error in the global DDM algorithm. Such is not the case
when recycling is not used, as was shown through numerical experiments
in the previous section.

\quad\\
Let us quantify more precisely the convergence rate.  Because
$\mathfrak{R}$ is symetric, its norm $\vert \mathfrak{R}\vert$ equals
its spectral radius. As a consequence, to bound the convergence rate
provided by the previous lemma, we need to estimate the largest
eigenvalue of $\mathfrak{R}$. This can be done explicitely by
examining its caracteristic polynomial.
\begin{equation*}
  \begin{aligned}
    \mrm{det}(\lambda\Id - \mathfrak{R})
    & = (\rho + 2\epsilon-\lambda)(2\epsilon-\lambda) - 4\epsilon (1+\epsilon)\\
    & = \lambda^{2}-\lambda (\rho + 4\epsilon) -(2-\rho)2\epsilon
  \end{aligned}
\end{equation*}
Because $\rho<1$, we see that $(2-\rho)2\epsilon>0$ and that the roots
of this polynomial have opposite signs. Hence $\vert\mathfrak{R}\vert$
agrees with the positive root.  With the gross estimate
$\sqrt{x+y}\leq \sqrt{x}+\sqrt{y}$, we deduce
\begin{equation}
  \begin{aligned}
    \vert \mathfrak{R}\vert
    & = \frac{1}{2}(\rho + 4\epsilon) +
    \frac{1}{2}\sqrt{(\rho + 4\epsilon)^2 +
      8(2-\rho)\epsilon},\\
    \vert \mathfrak{R}\vert
    & \leq  \rho+4\epsilon + 2\sqrt{\epsilon}.
  \end{aligned}
\end{equation}
To simplify the above estimate observe that, if
$\rho+4\sqrt{\epsilon}<1$, then $2\sqrt{\epsilon}<1\Rightarrow
4\epsilon < 2\sqrt{\epsilon}$ and in this case $\vert
\mathfrak{R}\vert\leq \rho+4\epsilon + 2\sqrt{\epsilon} <
\rho+4\sqrt{\epsilon}<1$. Let us examine what does the condition
$\rho+4\sqrt{\epsilon}<1$ means.  Recall that the operator
$(1-\alpha)\Id - \alpha\Pi\mS$ was proved to be a contraction, with
the estimate $\rho \leq 1-\alpha(1-\alpha)\gamma_h^2$, see
\cite[Thm.9.2]{claeys2020robust}. As a consequence, to ensure
$\rho+4\sqrt{\epsilon}<1$, it is sufficient that 
\begin{equation}\label{ConditionOnk}
  \epsilon = 2\Big(\frac{ \sqrt{\mrm{cond}(\Prec\mR^*\mT\mR)}-1}{
    (\sqrt{\mrm{cond}(\Prec\mR^*\mT\mR)}+1}\Big)^k <
  \big(\alpha(1-\alpha)\gamma_h^{2}/4\big)^2
\end{equation}
Because $\epsilon$ decays exponentially fast to $0$ as $k\to \infty$,
which reflects the spectral convergence of the (preconditioned)
conjugate gradient, only a few PCG iterations are necessary for
satisfying \eqref{ConditionOnk}. This is particularly true when the
preconditioner $\Prec$ is devised appropriately. The next lemma
summarizes the previous dsicussion on convergence criterion and
convergence rate.

\begin{lem}\quad\\
  Assume the number $k$ of PCG iterations constant and chosen
  sufficiently large to satisfy Condition \eqref{ConditionOnk}.  Then
  $\rho+4\sqrt{\epsilon}<1$, and the sequence $\tilde{\bq}^{(n)}$
  defined by \eqref{InexactRichardson} converges toward
  $\bq^{(\infty)}$ solution to \eqref{SkeletonFormulation} with
  the error estimate
  \begin{equation*}
    \Vert \tilde{\bq}^{(n)}-\bq^{(\infty)}\Vert_{\mT^{-1}}\leq
    \frac{\Vert \tilde{\bq}^{(0)}\Vert_{\mT^{-1}}}{1-(\rho+4\sqrt{\epsilon})}
    \;(\rho+4\sqrt{\epsilon})^{n}.
  \end{equation*}
\end{lem}

%%%%%%%%%%%%%%%%%%%%%%%%%%%%%%%%%%%%%%
%% \bibliography{draft}
%% \bibliographystyle{plain}
%%%%%%%%%%%%%%%%%%%%%%%%%%%%%%%%%%%%%%

\end{document}